\documentclass{article}
\usepackage{amsthm}
\usepackage{amssymb}
\usepackage{amsmath}
\usepackage{amsfonts}
\usepackage{color}
\usepackage{indentfirst} 
\usepackage{graphicx}
\usepackage[square, comma, sort&compress, numbers]{natbib}

\newtheorem{remark}{Remark}[section]

\newtheorem{theorem}{Theorem}[section]

\newtheorem{corollary}{Corollary}[section]

\def\b1{\mbox{\boldmath $1$}}

\newenvironment{demo*}{\vspace{3mm}\noindent{\bf Proof.}}{\hfill $\Box$ \vspace{3mm}}

\begin{document}
\title{\bf \Large {A New Class of Symmetric  Distributions   Including  the   Elliptically Symmetric Logistic}}
{\color{red}{\author{
\normalsize{Chuancun Yin\;\; Xiuyan Sha}\\
{\normalsize\it  School of Statistics,  Qufu Normal University}\\
\noindent{\normalsize\it Shandong 273165, China}\\
e-mail:  ccyin@qfnu.edu.cn}}}
\maketitle
\vskip0.01cm
\centerline{\large {\bf Abstract}}   We introduce a  new broad and flexible class of multivariate  elliptically symmetric distributions  including  the   elliptically symmetric logistic and multivariate normal.  Various probabilistic
properties of the new distribution are studied, including  the distribution of linear transformations, marginal distributions, conditional distributions,  moments, stochastic representations
and  characteristic function.

\medskip

\numberwithin{equation}{section}
\section{Introduction}\label{intro}
The elliptically symmetric logistic distribution with density
\begin{equation*}
f({\bf x})=\frac{\Gamma(\frac{n}{2})|\Sigma|^{-\frac{1}{2}}}{\pi^{\frac{n}{2}}\int_0^{\infty}u^{\frac{n}{2}-1} \frac{\exp(-u)}{(1+\exp(-u))^2}du}\frac{\exp(- ({\bf x}- {\boldsymbol \mu})^{T}{\bf \Sigma}^{-1}({\bf x-{\boldsymbol \mu}}))}{(1+\exp(- ({\bf x}- {\boldsymbol \mu})^{T}{\bf \Sigma}^{-1}({\bf x-{\boldsymbol \mu}})))^{2}},\; {\bf x}\in \Bbb{R}^n,
\end{equation*}
    has been originally introduced by Jensen (1985)  (as the generalization
of the multivariate normal distribution) and studied   by  Fang et al. (1990), Arnold (1992), Kano (1994),  Volodin (1999) and  G\'omez-S\'anchez-Manzano et al. (2006). Several applications  of  multivariate symmetric logistic distribution in risk management,  quantitative finance and actuarial science can be found in many literatures such as   Valdez and  Chernih(2003), Landsman (2004),  Landsman and Valdez (2003) and  Landsman et al. (2003, 2016, 2018), Xiao and Valdez (2015).  Note that several authors (cf. Gumbel (1961) Malik and Abraham (1973),   Ali et al. (1978), Fang
and Xu (1989),  Kotz, Balakrishnan and  Johnson (2000), Yeh (2010),  Hu and Lin (2018), Ghosh and Alzaatreh (2018)) have studied the multivariate logistic distribution using different definitions.

The elliptically symmetric logistic distribution   belongs to the elliptically contoured distributions family  (also called an elliptically symmetric
distributions family) $Ell_n({\boldsymbol \mu},{\bf \Sigma},g)$  with the location   parameter ${\boldsymbol \mu}$,   the  scale parameter  ${\bf \Sigma}$  and the density generator
$$g(u)=\frac{\exp(-u)}{(1+\exp(-u))^2}.$$
However, research work on the   multivariate symmetric logistic distribution    in recent decade is rather scarce    compared to much  research has focused on other  elliptical  distributions such as multivariate normal, multivariate Student $t$, multivariate Cauchy, multivariate power exponential distribution,  Kotz-type distribution and multivariate skew-normal distributions; see the books and papers of  Fang, Kotz and Ng (1990), Johnson, Kotz and Balakrishnan (1995) and Kotz, Balakrishnan and Johnson (2000), Azzalini and Regoli (2002),  Arellano-Valle and Azzalini (2013),  Battey and Linton (2014),  Arashi and Nadarajah (2016), Arellano-Valle et al.(2018) and the references therein.
    The relative few literatures on  the properties of elliptically symmetric logistic distribution  drive the authors to further study   the
properties of this distribution. Moreover,  we will define a new family of multivariate distributions   including the   elliptically symmetric logistic and multivariate normal and study the probabilistic properties of the distributions included by this family.


In some  literature, the   joint pdf of $n$-dimensional   elliptically symmetric logistic distribution  is defined as
\begin{equation*}
f({\bf x})=c_n|\Sigma|^{-\frac{1}{2}}\frac{\exp(-\frac12q({\bf x}))}{(1+\exp(-\frac12q({\bf x})))^{2}}, {\bf x}\in \Bbb{R}^n,
\end{equation*}
where $q({\bf x})=({\bf x}- {\boldsymbol \mu})^{T}{\bf \Sigma}^{-1}({\bf x-{\boldsymbol \mu}})$.
The density generator is
$$g_n(u)=\frac{\exp(-u)}{(1+\exp(-u))^2},$$
and the normalizing constant $c_n$ is given by
\begin{equation}
c_n=\frac{\Gamma(\frac{n}{2})}{(2\pi)^{\frac{n}{2}}}\left[\int_0^{\infty}x^{\frac{n}{2}-1}g_n(x)dx\right]^{-1}.
\end{equation}
Interested readers may refer to   Landsman and Valdez (2003) and  Landsman et al. (2003, 2016, 2018), for more details on the   elliptically symmetric logistic distribution and its applications.
As pointed out by Landsman and Valdez (2003) this normalizing constant has been mistakenly pointed in Fang et al. (1990), Gupta et al. (2013), Xiao and Valdez (2015).
Further simplification of the normalizing constant $c_n$ suggests by Landsman and Valdez (2003):
\begin{equation}
c_n=(2\pi)^{-\frac{n}{2}}\left[\sum_{j=1}^{\infty}(-1)^{j-1}j^{1-\frac{n}{2}}\right]^{-1}
\end{equation}
by using the expansion
\begin{equation}
\frac{e^{-x}}{(1+e^{-x})^2}=\sum_{j=1}^{\infty}(-1)^{j-1}je^{-j x},\; x>0.
\end{equation}
We observe that the formula (1.2) has no meaning when $n=1$ and 2, since the series
$$\sum_{j=1}^{\infty}(-1)^{j-1}\sqrt{j}\;\; {\rm and}\;\; \sum_{j=1}^{\infty}(-1)^{j-1}$$
are divergent.

We now give the definition of a generalized elliptically symmetric logistic
distribution  including the elliptically symmetric logistic distribution and multivariate normal.

{\bf Definition 1.1} The $n$-dimensional random vector $\bf{X}$ is said to have a  generalized elliptically symmetric logistic
distribution with location parameter  $\boldsymbol{\mu}$ ($n$-dimensional vector) and dispersion matrix $\bf{\Sigma}$ ($n\times n$ matrix with $\bf{\Sigma}>0$) if
    its pdf has the form
\begin{equation}
f({\bf x})=d_n|{\bf \Sigma}|^{-\frac{1}{2}}g(({\bf x}- {\boldsymbol \mu})^{T}{\bf \Sigma}^{-1}({\bf x-{\boldsymbol \mu} })), \;{\bf x}\in \Bbb{R}^n,
\end{equation}
 where $d_n$ is the normalizing constant and will be determined in the next section,  $a,b,r >0$ are constants and
 \begin{equation}
 g(u)=\frac{\exp(-bu)}{(1+\exp(-au))^r},
 \end{equation}
 is its  density generator. If $\bf X$ belongs to the generalized elliptically symmetric logistic distribution, we shall write ${\bf{X}}\sim GML_n ({\boldsymbol \mu},{\bf \Sigma},g)$.

The  generalized elliptically symmetric logistic distribution  is a particular case of an elliptically distribution, so $\bf X$ admits the stochastic representation
\begin{equation}
{\bf X}={\boldsymbol \mu}+\sqrt{R}{\bf A}'{\bf U}^{(n)},
\end{equation}
where ${\bf A}$  is a square matrix such that ${\bf A}'{\bf A}= {\bf \Sigma}$, ${\bf U}^{(n)}$ is uniformly distributed on the unit sphere surface in $\Bbb{R}^n$, and $R\ge 0$ is independent of ${\bf U}^{(n)}$
and has the pdf given by
\begin{equation}
f_R(v)=\frac{1}{\int_0^{\infty}t^{\frac{n}{2}-1}g(t)dt}v^{\frac{n}{2}-1}g(v), v\ge 0.
\end{equation}
For more details see Cambanis et al. (1981).  Note that the  $n$-dimensional   elliptically symmetric logistic distribution    can be deduced as a special case of (1.4) by setting $a=b=1$ and $r=2$; For $b=\frac12$ and $r=0$,   we obtain the multivariate normal distribution;
The  density generator (1.5) covers  all the density generators of the generalized logistic types I-IV distributions in Arashi and Nadarajah (2017).
Figures 1-5 illustrate the  density functions of (1.4) with $n=2,a=b=1$ and  $r=0.5,1,2,5,10$.
\begin{figure}[htbp]
\centering
\begin{minipage}[t]{0.9\textwidth}
\centering
\includegraphics[width=10cm]{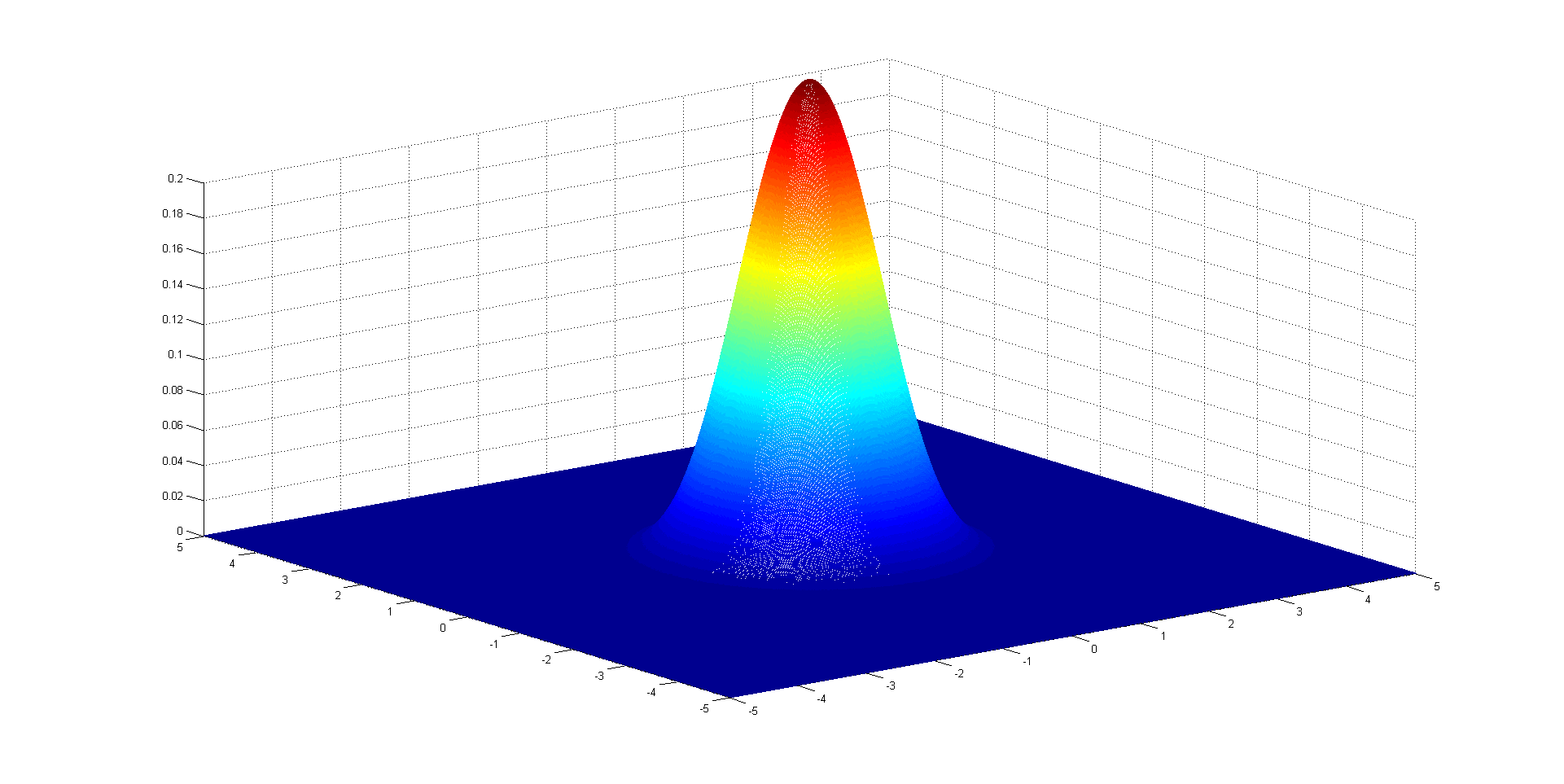}
\caption{ $r=0.5$}
\end{minipage}
\begin{minipage}[t]{0.9\textwidth}
\centering
\includegraphics[width=10cm]{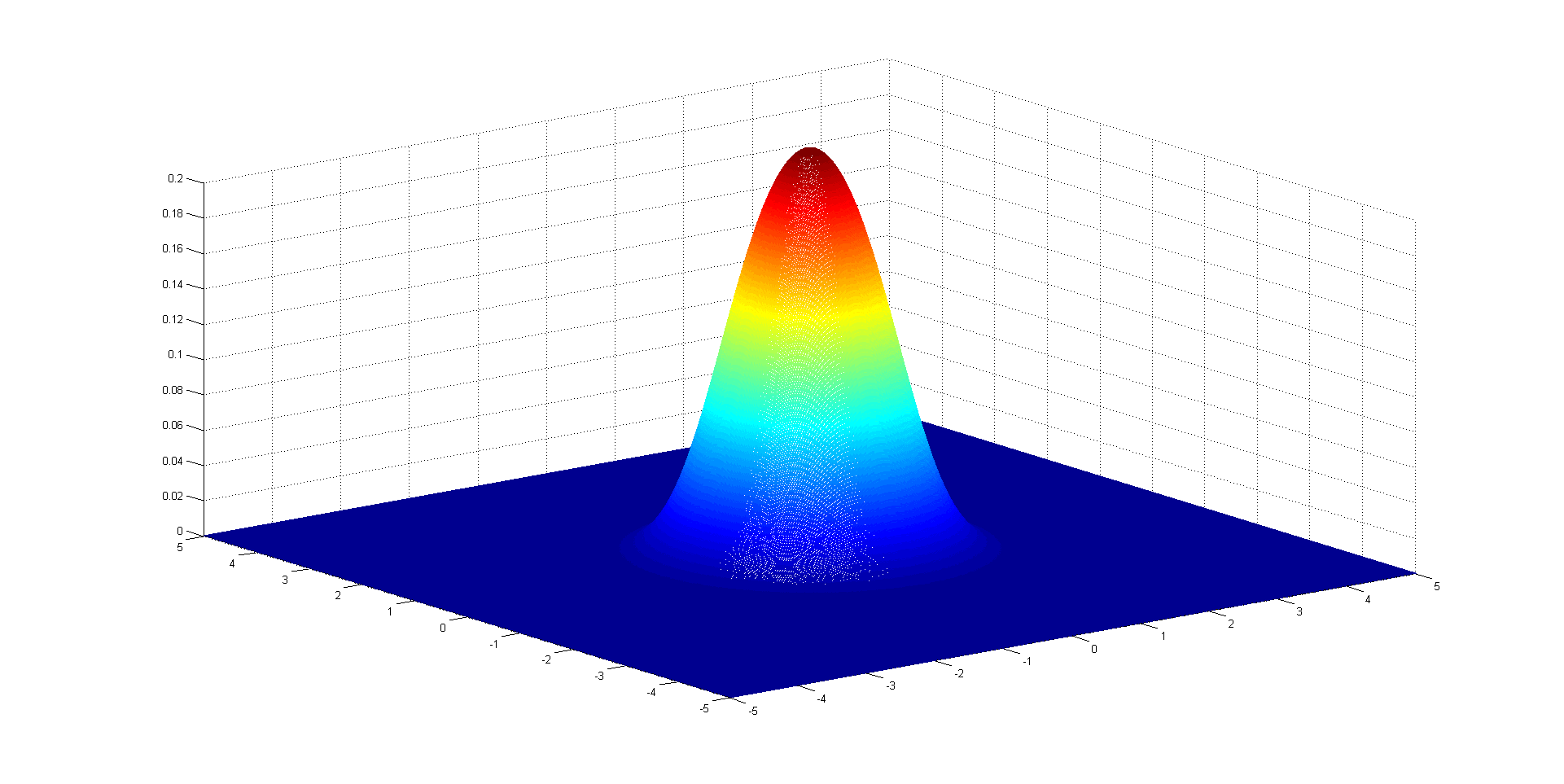}
\caption{$r=1$}
\end{minipage}
\begin{minipage}[t]{0.9\textwidth}
\centering
\includegraphics[width=10cm]{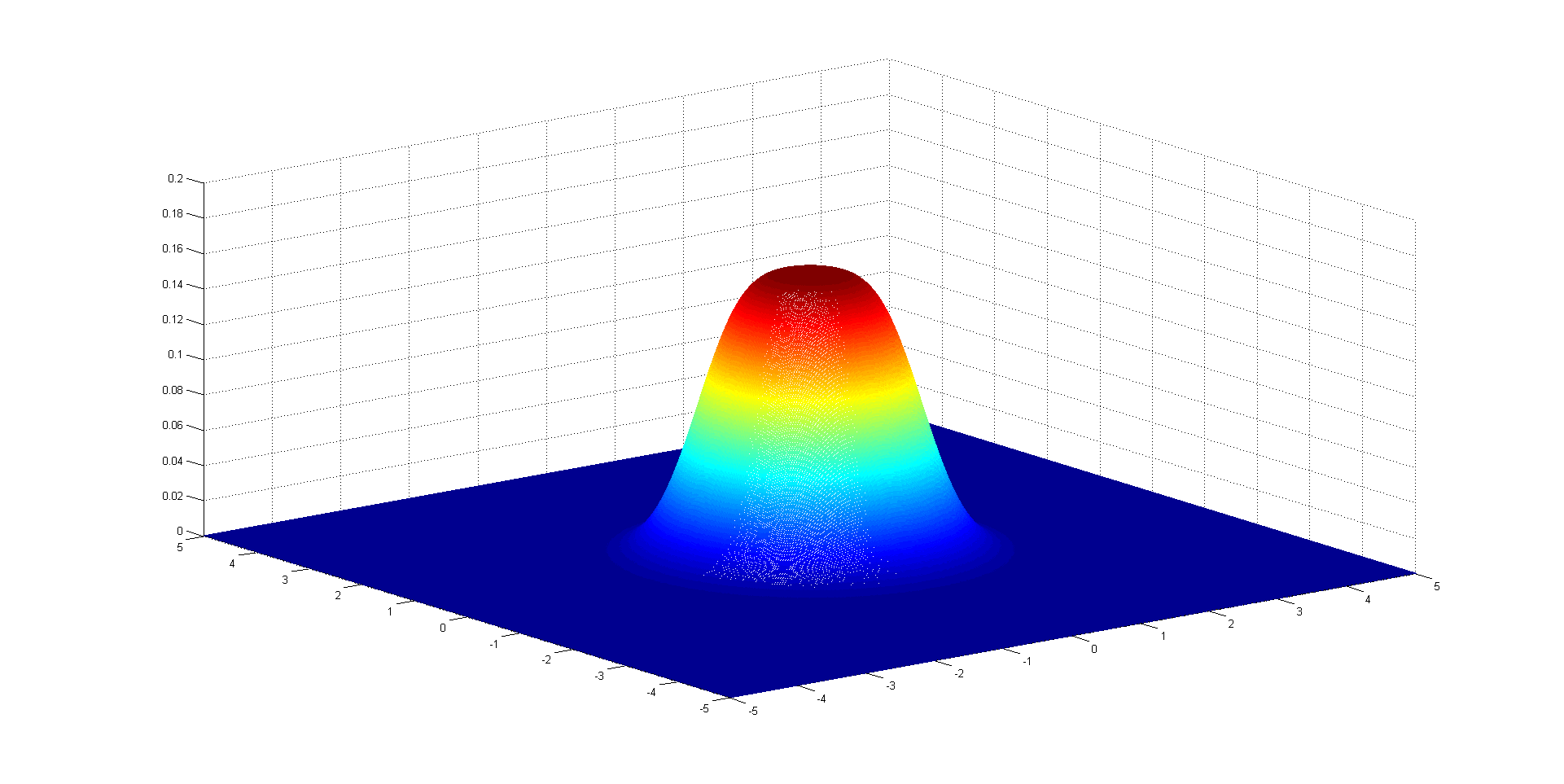}
\caption{$r=2$}
\end{minipage}
\begin{minipage}[t]{0.9\textwidth}
\centering
\includegraphics[width=10cm]{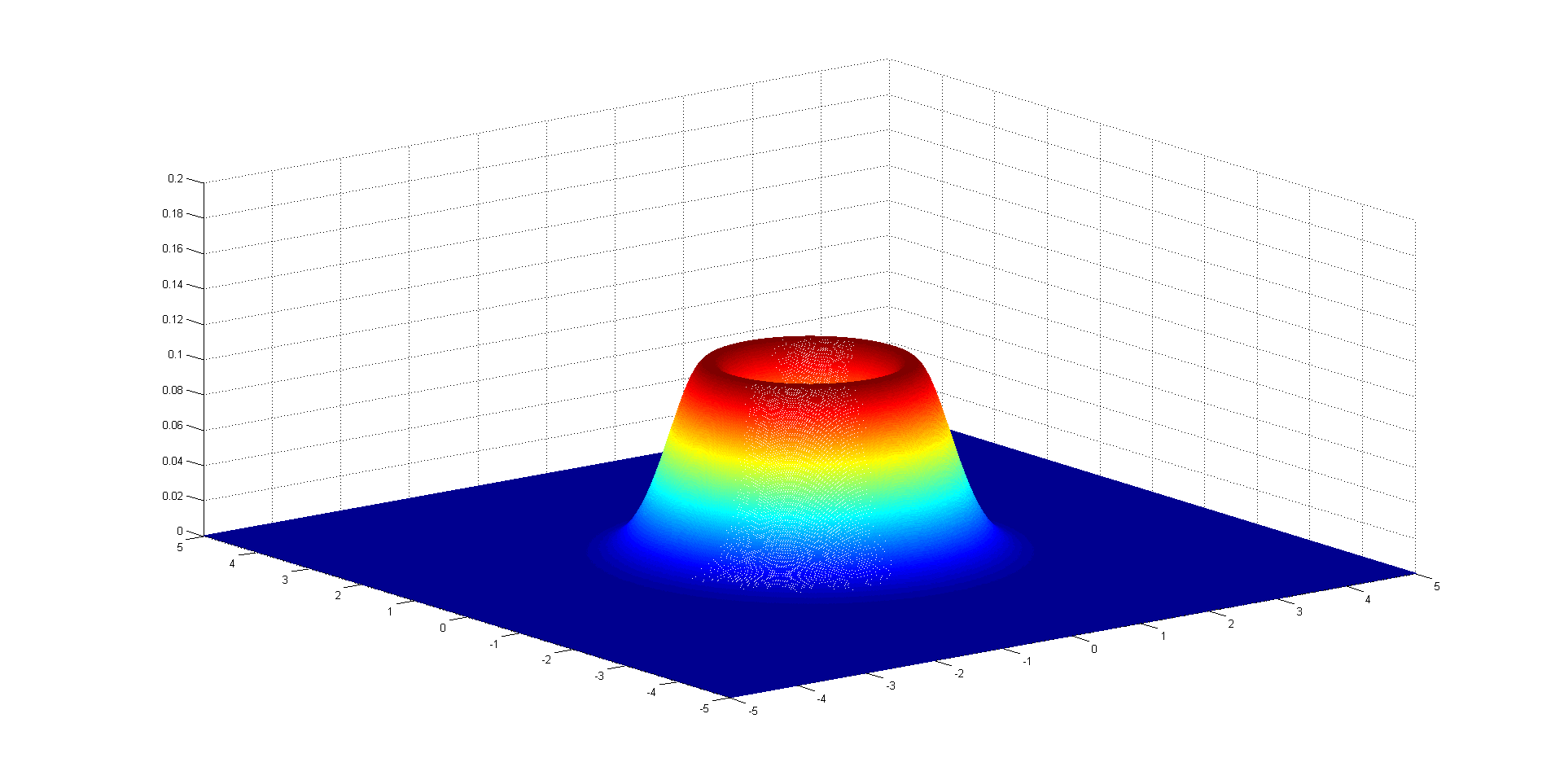}
\caption{$r=5$}
\end{minipage}
\end{figure}
\begin{figure}
\begin{minipage}[t]{0.9\textwidth}
\centering
\includegraphics[width=10cm]{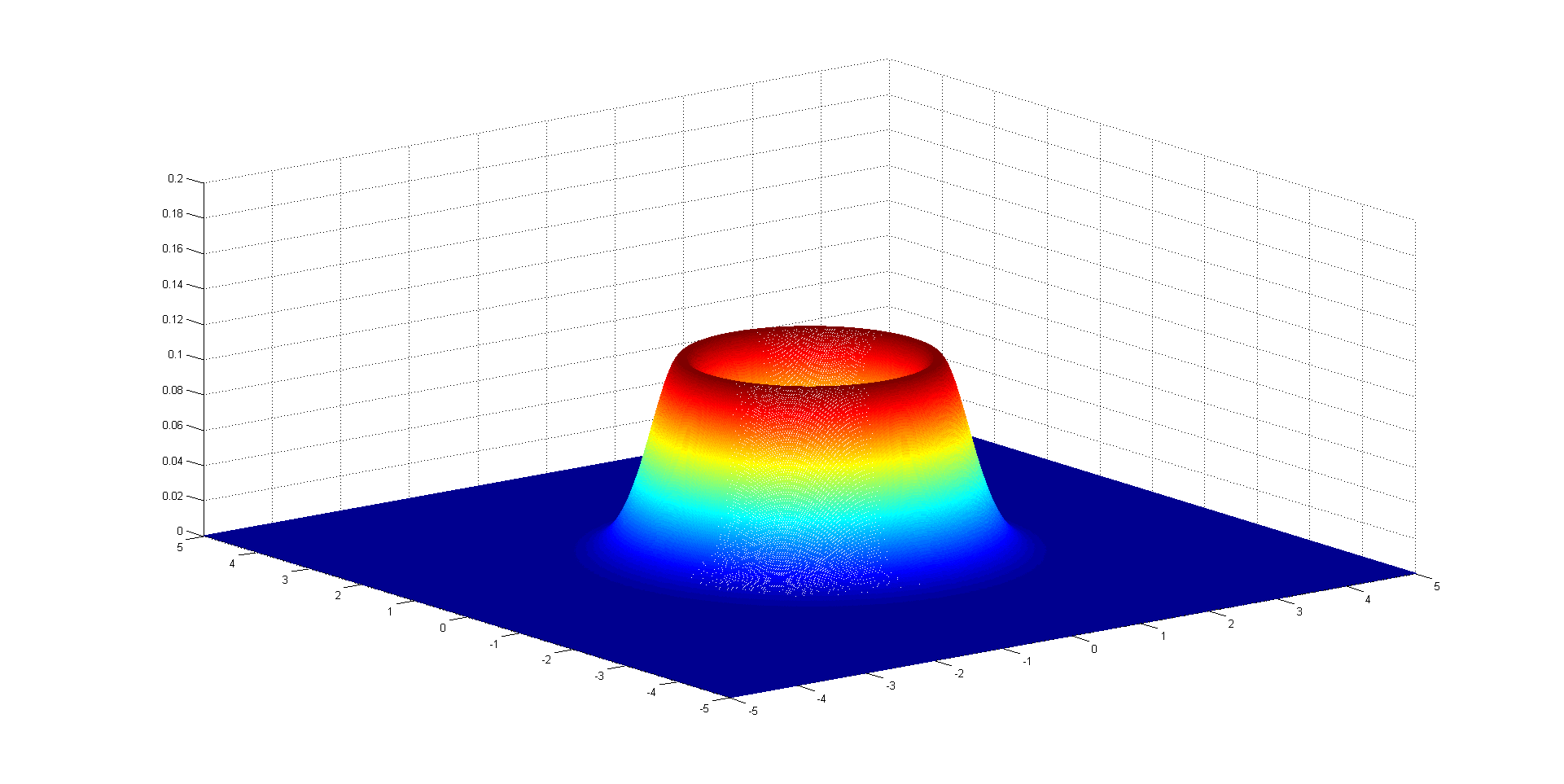}
\caption{$r=10$}
\end{minipage}
\end{figure}

The rest of the paper is organized as follows. In Section 2, we discuss the expression of the  normalizing constant in (1.1) and   give the correct values for $c_1$ and $c_2$, and provide a  equivalent expression for $c_n, n\ge 3.$  In Sections 3-6 we will investigate some of the properties  of this new class of the multivariate distributions.  More specifically, we find   the distributions of   linear transformations, the  marginal and conditional distributions, the  moments and the characteristic functions.

 \section{Evaluation of  the  normalizing constants }

We first collect some facts  of   the Riemann zeta function and  the generalized Hurwitz-Lerch Zeta function which will be used in the sequel.
The Riemann Zeta function $\zeta$ is defined as
\begin{eqnarray*}
 \zeta(s)=\left\{\begin{array}{ll} \sum_{n=1}^{\infty}\frac{1}{n^s}=\frac{1}{1-2^{-s}}\sum_{n=1}^{\infty}\frac{1}{(2n-1)^s},  \ & {\rm if}\; {\cal{R}}(s)>1,\\
\frac{1}{1-2^{1-s}}\sum_{n=1}^{\infty}\frac{(-1)^{n+1}}{n^s} ,\ &{\rm if}\;   {\cal{R}}(s)>0, s\neq 1,
 \end{array}
  \right.
\end{eqnarray*}
which can, except for a simple pole at $s=1$ with its residue 1, be continued
meromorphically to the whole complex $s$-plane; see, for details  Srivastava(2003),  Choi et al. (2004).
 Recall that  (see Arakawa et al. (2014), Cvijovi\'c and Klinowski (2002))
\begin{equation}
\zeta(2n)=(-1)^{n+1}\frac{2^{2n-1}}{(2n)!}\pi^{2n}B_{2n}, n\in\Bbb{N}_0,
\end{equation}
and
\begin{equation}
\zeta(2n+1)=(-1)^{n+1}\frac{(2\pi)^{2n+1}}{2(2n+1)!}\int_0^1 B_{2n+1}(u)\cot(\pi u)du, n\in\Bbb{N},
\end{equation}
where  $B_n=B_n(0)$ are  the   $n$th Bernoulli numbers and  $B_n(x)$   are  Bernoulli   polynomials  defined by the generating function
$$\frac{te^{tx}}{e^t-1}=\sum_{n=0}^{\infty}B_n(x)\frac{t^n}{n!}, \; |t|<2\pi.$$
The Bernoulli numbers are well-tabulated (see, for example, Srivastava(2003)):
$$B_0=1,B_1=-\frac12, B_2=\frac16, B_4=-\frac{1}{30}, B_6=\frac{1}{42},B_{2n+1}=0 \;(n=1,2,\cdots),\cdots.$$

The functions $\zeta(s)$   has the following integral
representations (cf.   Srivastava and  Choi (2012, p.169, p.172))
$$\zeta(s)=\frac{1}{\Gamma(s)}\int_0^{\infty}\frac{t^{s-1}}{e^t-1}dt,\;\; {\cal{R}}(s)>1,$$
and
$$\zeta(s)=\frac{(1-2^{1-s})^{-1}}{\Gamma(s+1)}\int_0^{\infty}\frac{t^{s}e^{t}}{(e^t+1)^2}dt,\;\; {\cal{R}}(s)>0.$$
Note that  there is an extra 2 in (51) of  Srivastava and  Choi (2012, p.172).

The generalized Hurwitz-Lerch Zeta function is defined by (cf. Lin et al.(2006))
$$\Phi^*_{\mu}(z,s,a)=\frac{1}{\Gamma(\mu)}\sum_{n=0}^{\infty}\frac{\Gamma(\mu+n)}{n!}\frac{z^n}{(n+a)^s},$$
which has an integral representation
\begin{equation}
\Phi^*_{\mu}(z,s,a)=\frac{1}{\Gamma(s)}\int_0^{\infty}\frac{t^{s-1}e^{-at}}{(1-ze^{-t})^{\mu}}dt,
\end{equation}
where ${\cal{R}}(a)>0$; ${\cal{R}}(s)>0$ when $|z|\le 1 (z\neq 1); {\cal{R}}(s)>1$ when $z=1).$

\begin{theorem} Consider the normalizing constant $d_n$ defined in (1.4).\\
 Then
\begin{equation}
d_n=\left(\frac{a}{\pi}\right)^{\frac{n}{2}}\left[\Phi^*_{r}\left(-1,\frac{n}{2},\frac{b}{a}\right)\right]^{-1},
\end{equation}
where $\Phi^*_{r}$ is the generalized Hurwitz-Lerch Zeta function.
\end{theorem}
{\bf Proof}\;  By using the formula (1.37) in Denuit et al. (2005), we have
\begin{eqnarray*}
 d_n&=&\frac{\Gamma(\frac{n}{2})}{(\pi)^{\frac{n}{2}}}\left[\int_0^{\infty}x^{\frac{n}{2}-1}g(x)dx\right]^{-1}\\
&=&\frac{\Gamma(\frac{n}{2})}{(\pi)^{\frac{n}{2}}}\left[\int_0^{\infty}x^{\frac{n}{2}-1}\frac{\exp(-bx)}{(1+\exp(-ax))^r}dx\right]^{-1}\\
&=&\frac{\Gamma(\frac{n}{2})a^{\frac{n}{2}}}{\pi^{\frac{n}{2}}}
\left[\int_0^{\infty}x^{\frac{n}{2}-1}\frac{\exp(-\frac{b}{a}x)}{(1+\exp(-x))^r}dx\right]^{-1}\\
&=& \left(\frac{a}{\pi}\right)^{\frac{n}{2}}\left[\Phi^*_{r}\left(-1,\frac{n}{2},\frac{b}{a}\right)\right]^{-1},
\end{eqnarray*}
as desired.
\begin{corollary} Consider the normalizing constant $c_n$ defined in (1.1).\\
(i)\ If $n=1$, then
$$c_1=\frac{1}{\sqrt{2\pi}}  \left( \Phi^*_{2}(-1,\frac12,1)\right)^{-1},$$
where $\Phi^*_{\mu}$ is the generalized Hurwitz-Lerch Zeta function.\\
(ii)\ If $n=2$, then $c_2=\pi^{-1}.$\\
(iii)\ If $n=4$, then $c_4=\frac{1}{4\pi^2\ln 2}.$\\
(iv)\ If $n\ge 3, n\neq 4$, then
$$c_n=\pi^{-\frac{n}{2}}\left[\left(2^{\frac{n}{2}}-4\right)\zeta\left(\frac{n}{2}-1\right)\right]^{-1},$$
where
$\zeta$ is the Riemann zeta function.
\end{corollary}
{\bf Proof}\; (i) The result follows by letting  $n=1, a=b=1, r=2$ in (2.4) and some algebras.\\
(ii) If $n=2$, by (1.1) we have
\begin{equation*}
c_n=\frac{\Gamma(\frac{2}{2})}{(2\pi)^{\frac{2}{2}}}\left[\int_0^{\infty}x^{\frac{2}{2}-1}\frac{\exp(-x)}{(1+\exp(-x))^2}dx\right]^{-1}=\frac{1}{\pi},
\end{equation*}
where we have used the fact that
$$h(x)=\frac{2\exp(-x)}{(1+\exp(-x))^2}, \; x>0,$$
is the pdf of the half-logistic distribution.\\
(iii) By (1.2),
$$c_4=(2\pi)^{-2}\left[\sum_{j=1}^{\infty}(-1)^{j-1}\frac{1}{j}\right]^{-1}= \frac{1}{4\pi^2\ln 2},$$
where we have used the well known fact
$$\sum_{n=1}^{\infty}(-1)^{n-1}\frac{1}{n}=\ln 2.$$
(iv) If $n\ge 3, n\neq 4$, by making use of
$$\zeta(s)=\frac{1}{1-2^{1-s}}\Phi^*_{2}(-1,s+1,1),\;  {\cal{R}}(s)>0,$$
we have
 \begin{eqnarray*}
c_n=(2\pi)^{-\frac{n}{2}}\left[\Phi^*_{2}(-1,\frac{n}{2},1)\right]^{-1}
=\pi^{-\frac{n}{2}}\left[\left(2^{\frac{n}{2}}-4\right)\zeta\left(\frac{n}{2}-1\right)\right]^{-1},
 \end{eqnarray*}
 where $\zeta$ is the Riemann zeta function.

\begin{remark} When $\frac{n}{2}-1$ are even or odd  positive integers, we apply formulae (2.1) and (2.2) to (2.4):
For $n=4m+2, m=1,2,\cdots$, the following formula holds,
$$c_n=\frac{(2m)!}{(-1)^{m+1}2^{2m-2}\pi^{\frac{n}{2}+2m}(2^{\frac{n}{2}}-4)B_{2m}};$$
For $n=4m+4, m=1,2,\cdots$, the following formula holds,
$$c_n=\frac{(2m+1)!}{2^{\frac{n}{2}-1} (2^{\frac{n}{2}}-4)(-1)^{m+1}(2\pi)^{2m+1}\int_0^1B_{2m+1}(u)\cot(\pi u)du}.$$
\end{remark}
\noindent {\bf Example 2.1.}\;  Let us compute $c_n$ for small  even $n$.
 By using the fact that (cf. Srivastava and  Choi (2012, p.167))
 $$\zeta(2)=\frac{\pi^2}{6}, \;\zeta(4)=\frac{\pi^4}{90},\;\zeta(6)=\frac{\pi^6}{945},\;\zeta(8)=\frac{\pi^8}{9450},$$
 we get
 $$c_6=\frac{ 3}{2}\pi^{-5},\; c_{10}=\frac{45}{14}\pi^{-9},\; c_{14}=\frac{945}{124}\pi^{-13}, \; c_{18}=\frac{4725}{254}\pi^{-17}.$$
\begin{remark} Using the relationship of $c_n$ and $d_n$ above we find the following formula
$$\Phi^*_{2}(-1,1,1)=\frac12,\; \; \Phi^*_{2}(-1,2,1)=\ln 2,$$
and
$$\Phi^*_{2}\left(-1,\frac{n}{2},1\right)=\left(2^{\frac{n}{2}}-4\right)\zeta\left(\frac{n}{2}-1\right), n\ge 3, n\neq 4.$$
\end{remark}

\section{Linear transformations}

Consider the affine transformations of the form ${\bf Y = BX+b}$ of a random vector  ${\bf{X}}\sim GML_n ({\boldsymbol \mu},{\bf \Sigma},g)$. If
${\bf B}$ is a nonsingular $n\times n$ matrix it can
be easily verified by definition or by the  characteristic function in Section 6 that   ${\bf{Y}}\sim GML_n ({\bf B}{\boldsymbol \mu}+{\bf b},{\bf B}{\bf \Sigma}{\bf B}',g)$. However,  when $\bf B$ is a $m\times n$ matrix with $m<n$ and $rank({\bf B})=m$,  the following theorem shows that the density generator of $\bf Y$ is not necessarily $g$, it may dependent on $n$ and $m$.
\begin{theorem} Let ${\bf{X}}\sim GML_n ({\boldsymbol \mu},{\bf \Sigma},g)$,  with stochastic representation ${\bf X}={\boldsymbol \mu}+\sqrt{R}{\bf A}'{\bf U}^{(n)}$. Let ${\bf Y = BX+b}$, where  $\bf B$ is a $m\times n$ matrix with $m<n$ and $rank({\bf B})=m$ and ${\bf b}\in\Bbb{R}^m$. Then  ${\bf Y} \sim Ell_n ({\bf B}{\boldsymbol \mu}+{\bf b}, {\bf B}{\bf \Sigma}{\bf B}',g_y)$ with
$$g_y(t)= e^{-bt}a^{-\frac{n-m}{2}}\Gamma\left(\frac{n-m}{2}\right)\Phi^*_{r}\left(-e^{-at},\frac{n-m}{2},\frac{b}{a}\right),$$
where $\Phi^*_{r}$ is the generalized Hurwitz-Lerch Zeta function.
Moreover, ${\bf Y}$ admits the stochastic representation
\begin{equation*}
{\bf Y}={\boldsymbol \mu}_Y+\sqrt{R}_Y{\bf A}_Y'{\bf U}^{(m)},
\end{equation*}
where ${\bf A}_Y'{\bf A}_Y={\bf B}{\bf \Sigma}{\bf B}'$ and $R_Y=R_X B$. Here $B$ is the nonnegative random variable independent of $R$ with distribution $Beta(\frac{m}{2},\frac{n-m}{2})$.
\end{theorem}
{\bf Proof} The result  follows from  Theorems 2.15 and 2.16
in Fang et al. (1990). In the sequel, we provided an  alternative proof. In dong so we consider the transformation
$${\bf Y^*}=\left(\begin{array}{c}
{\bf Y_1}\\
{\bf Y_2}
\end{array}\right)=\left(\begin{array}{c}
{\bf B}\\
{\bf C}
\end{array}\right){\bf X}+\left(\begin{array}{c}
{\bf b}\\
{\bf 0}_{n-m}
\end{array}\right),$$
where ${\bf C}$ is any given matrix such that $\left(\begin{array}{c}
{\bf B}\\
{\bf C}
\end{array}\right)$ is nonsingular. So that
$${\bf{Y}^*}\sim GML_n \left(\left(\begin{array}{c}
 {\bf B}{\boldsymbol \mu}+{\bf b}\\
{\bf C{\boldsymbol \mu}}
\end{array}\right),
 \left(\begin{array}{cc}
 {\bf B}{\bf \Sigma}{\bf B}' \; {\bf B}{\bf \Sigma}{\bf C}'\\
 {\bf C}{\bf \Sigma}{\bf B}'\; {\bf C}{\bf \Sigma}{\bf C}'\\
\end{array}\right),g\right),$$
The result follows, since
\begin{eqnarray*}
\int_{\Bbb{R}^{n-m}}g\left(\sum_{i=1}^n x_i^2\right)\prod_{j=m+1}^n dx_{j}&=&2^{n-m}\int_{\Bbb{R}_{+}^{n-m}}g\left(\sum_{i=1}^n x_i^2\right)dx_{m+1}\cdots dx_n\\
&=&\int_{\Bbb{R}_{+}^{n-m}}g\left(\sum_{i=1}^m x_i^2+\sum_{j=m+1}^n u_j\right)\prod_{j=m+1}^n u_j^{-\frac12}du_{m+1}\cdots du_n\\
&=&\int_{D_1}g\left(\sum_{i=1}^m x_i^2+y_n\right)\prod_{j=m+1}^{n-1} y_j^{-\frac12}\left(y_n-\sum_{j=m+1}^{n-1}y_j\right)^{-\frac12}dy_{m+1}\cdots dy_n\\
&=&\int_{D_2}\left(1-\sum_{j=m+1}^{n-1}v_j\right)^{-\frac12}\prod_{j=m+1}^{n-1}\left(v_j^{-\frac12}dv_j\right)
\int_0^{\infty}g\left(\sum_{i=1}^m x_i^2+y\right)y^{\frac{n-m}{2}-1}dy\\
&=&\frac{\pi^\frac{n-m}{2}}{\Gamma(\frac{n-m}{2})}\int_0^{\infty}g\left(\sum_{i=1}^m x_i^2+y\right)y^{\frac{n-m}{2}-1}dy,
\end{eqnarray*}
where
$$D_1=\{(y_{m+1},\cdots, y_n)|y_i\ge 0,i=m+1,\cdots,n, \sum_{j=m+1}^{n-1}y_i\le y_n\},$$
and
$$D_2=\{(v_{m+1},\cdots, v_{n-1})|v_i\ge 0,i=m+1,\cdots,n-1, \sum_{j=m+1}^{n-1}v_i\le 1\}.$$
This ends the proof.

Taking $B=(\alpha_1,\cdots,\alpha_n):={\boldsymbol \alpha}'$ in Theorem 3.1 leads to
$${\boldsymbol \alpha}'{\bf X}\sim Ell_1({\boldsymbol \alpha}'{\boldsymbol \mu,  {\boldsymbol \alpha}'{\bf \Sigma}{\boldsymbol \alpha},g_{1,n}}).$$
In particular,
$$X_k\sim Ell_1(\mu_k,\sigma_k^2, g_{1,n}), \;\; k=1,2,\cdots, n,$$
 and
 $$\sum_{k=1}^n X_k\sim Ell_1\left(\sum_{k=1}^n\mu_k,\sum_{k=1}^n\sum_{l=1}^n \sigma{_kl}, g_{1,n}\right),$$
where
$$g_{1,n}(t)=a^{-\frac{n-1}{2}}\Gamma\left(\frac{n-1}{2}\right)\Phi^*_{r}\left(-e^{-at},\frac{n-1}{2},\frac{b}{a}\right).$$

\section{ Marginal and conditional distributions}

For fixed $m<n$, consider the partitions of  $\bf X$, ${\boldsymbol \mu}, {\bf \Sigma}$  given below
$${\bf X}=\left(\begin{array}{c}
{\bf X^{(1)}}\\
{\bf X^{(2)}}
\end{array}\right), \;
\; {\boldsymbol \mu}=\left(\begin{array}{c}
{{\boldsymbol \mu}^{(1)}}\\
{{\boldsymbol \mu}^{(2)}}
\end{array}\right),\;
{\bf \Sigma}=\left(\begin{array}{cc}
 {\bf \Sigma}_{11} \; {\bf \Sigma}_{12}\\
 {\bf \Sigma}_{21} \; {\bf \Sigma}_{22}\\
\end{array}\right),
$$
where ${\bf X^{(1)}}, {{\boldsymbol \mu}^{(1)}}\in\Bbb{R}^m$ $(m<n)$, ${\bf X^{(2)}}, {{\boldsymbol \mu}^{(2)}}\in\Bbb{R}^{n-m}$,
${\bf \Sigma}_{11}$ is $m\times m$ matrix, ${\bf \Sigma}_{12}$ is $m\times (n-m)$ matrix,  ${\bf \Sigma}_{21}$ is $(n-m)\times m$ matrix,
${\bf \Sigma}_{22}$ is $(n-m)\times (n-m)$ matrix.

The following theorem gives the result on the marginal distributions of ${\bf X^{(1)}}$ and ${\bf X^{(2)}}$.

\begin{theorem} Let ${\bf{X}}\sim GML_n ({\boldsymbol \mu},{\bf \Sigma},g)$, where $g$ is defined as (1.6), then\\
(i)\; $${\bf{X}^{(1)}}\sim Ell_m ({\boldsymbol \mu}^{(1)},{\bf \Sigma}_{11},g_{(m)}),$$
(ii)\; $${\bf{X}^{(2)}}\sim Ell_{n-m} ({\boldsymbol \mu}^{(2)},{\bf \Sigma}_{22},g_{(n-m)}),$$
where $g_{(k)}$ is the function given by
$$g_{(k)}(t)= e^{-bt}a^{-\frac{n-k}{2}}\Gamma\left(\frac{n-k}{2}\right)\Phi^*_{r}\left(-e^{-at},\frac{n-k}{2},\frac{b}{a}\right),$$
where $\Phi^*_{r}$ is the generalized Hurwitz-Lerch Zeta function.

In particular, if $m = n - 2$ and $a=b$, then
\begin{eqnarray*}
 g_{(n-2)}(t)=\left\{\begin{array}{ll} \frac{1}{a(r-1)}\left(1-\frac{1}{(1+e^{-at})^{r-1}}\right),  \ & {\rm if}\; r\neq 1,\\
\frac{1}{a}\ln(1+e^{-at}),\ &{\rm if}\;    r=1.
 \end{array}
  \right.
\end{eqnarray*}
\end{theorem}
{\bf Proof}.\; Taking ${\bf B}=({\bf I}_m \;{\bf 0}_{m\times(n-m)})$ in Theorem 3.1, where ${\bf I}_m$ is the  $m\times m$ identity matrix and ${\bf 0}_{m\times(n-m)}$ is the $m\times(n-m)$ null matrix, we get ${\bf X}^{(1)}={\bf B}{\bf X}$. Applying Theorem  3.1 the result (i) follows. The proof of (ii) is similarly.

 An interesting question is to know whether the marginal pdfs of (1.4) have (1.5) as their density generator.
From Theorem 4.1 we see that  the marginal pdfs of (1.4) depending on the dimension $n$ and have not (1.5) as their density generator. Conversely, if an $n$-dimensional multivariate  elliptically symmetric  distribution with  $p$-dimensional marginal  generalized elliptically symmetric logistic, then the density generator $g$ is determined by
\begin{equation}
\frac{\exp(-bu)}{(1+\exp(-au))^r}=\int_t^{\infty}(w-t)^{\frac{n-p}{2}-1}g(w)dw.
\end{equation}
For example, if $a=b=1,n=2,r=2,p=1$, then  (4.1) becomes
\begin{equation}
\frac{\exp(-t)}{(1+\exp(-t))^2}=\int_0^{\infty}w^{-\frac{1}{2}}g(t+w)dw,
\end{equation}
from which we get
\begin{equation}
g(t)=\sum_{k=1}^{\infty}(-1)^{k-1}k^{\frac{3}{2}}e^{-kt}, \; t>0.
\end{equation}
The $2$-dimensional multivariate  elliptically symmetric  distribution with   the  density generator (4.3)  is not a  generalized elliptically symmetric logistic distribution.
If $a=b=1,n=3,r=2,p=1$, then  (4.1) becomes
\begin{equation}
\frac{\exp(-t)}{(1+\exp(-t))^2}=\int_t^{\infty}g(w)dw,
\end{equation}
 from which we get
\begin{equation}
g(t)=\frac{e^{-t}-e^{-2t}}{(1+e^{-t})^3}, \; t>0.
\end{equation}
 The $3$-dimensional multivariate  elliptically symmetric  distribution with   the  density generator (4.5)  is not a  generalized elliptically symmetric logistic distribution.

Therefore,  the  distribution (1.4) is not dimensionally coherent or consistent. A spherical distribution with density generator $f$ is said to have the consistency property if
\begin{equation}
\int_{-\infty}^{\infty}f\left(\sum_{i=1}^{n+1}x_i^2\right)dx_{n+1}=f\left(\sum_{i=1}^{n}x_i^2\right)
\end{equation}
for any integer $n$ and almost aa ${\bf x}\in\Bbb{R}^n$. This consistency property ensures that any marginal distribution of ${\bf X}$ also belongs to the same spherical family. Kano (1994) gave several necessary and sufficient conditions for a spherical distribution to satisfy (4.6) and list some examples. Dimensionally coherent elliptically distributions are multivariate  Normal, multivariate Student $t$, multivariate Cauchy and symmetric
stable. Distributions that do not have this property include the multivariate Logistic $\{\frac{e^{-u}}{1+e^{-u}}\}$, Pearson type II $\{(1-u)^{v-1}\}$ with $v>1$, Pearson type VII $\{(1+u)^{-v}\}$ with $v>\frac{n}{2}$, Kotz type $\{u^{N-1}e^{-r u^s}\}$ with $r,s>0, 2N+n>2$, and  multivariate Bessel distributions.

Voldin(1999) provided exact formulae for pdf of spherically symmetric distribution with logistic marginals. Applying  formulae (1) and (2) in Voldin(1999)  to   elliptically symmetric logistic distribution with
density generator (1.5) yields the density generators of  pdf of spherically symmetric distribution with  generalized elliptically  symmetric logistic marginals. That is
for $n=2m+1, m\in\Bbb{N}^+$,
\begin{equation}
g_n(t)=\frac{(-1)^m}{\pi^m}\frac{\partial ^m}{\partial t^m}\left(\frac{\exp(-bt)}{(1+\exp(-at))^r} \right), t\ge 0,
\end{equation}
and for $n=2m, m\in\Bbb{N}^+$,
\begin{equation}
g_n(t)=\frac{(-1)^m}{\pi^m}\frac{\partial ^m}{\partial t^m}\int_0^{\infty}\left(\frac{\exp(-b(t+z^2))}{(1+\exp(-a(t+z^2)))^r} \right)dz, t\ge 0.
\end{equation}
The following theorem gives the conditional distribution of ${\bf X}^{(2)}$ given  ${\bf X}^{(1)}$.
\begin{theorem} Let ${\bf{X}}\sim GML_n ({\boldsymbol \mu},{\bf \Sigma},g)$, where $g$ is defined as (1.5). Conditionally on ${\bf X}^{(1)}={\bf x}^{(1)}$, we have the conditional distribution of ${\bf X}^{(2)}$ is the elliptical distribution
$Ell_{n-m} ({\boldsymbol \mu}_{2.1},{\bf \Sigma}_{22.1}, g_{(2.1)})$,
where
$${\boldsymbol \mu}_{2.1}={\boldsymbol \mu}^{(2)}+{\bf \Sigma}_{21}{\bf \Sigma}_{11}^{-1}({\bf x}^{(1)}-{\boldsymbol \mu}^{(1)}),$$
$${\bf \Sigma}_{22.1}={\bf \Sigma}_{22}-{\bf \Sigma}_{21}{\bf \Sigma}_{11}^{-1}{\bf \Sigma}_{12},$$
and
$$ g_{(2.1)}(t)=\frac{e^{-b[t+({\bf x}^{(1)}-{\boldsymbol \mu}^{(1)})'{\bf \Sigma}_{11}^{-1}({\bf x}^{(1)}-{\boldsymbol \mu}^{(1)})]}}
{(1+e^{-a[t+({\bf x}^{(1)}-{\boldsymbol \mu}^{(1)})'{\bf \Sigma}_{11}^{-1}({\bf x}^{(1)}-{\boldsymbol \mu}^{(1)})]})^r},\; t\ge 0.$$
\end{theorem}
{\bf Proof}. The conditional density of  ${\bf X}^{(2)}$ given  ${\bf X}^{(1)}={\bf x}^{(1)}$ is given by
$$f({\bf x}^{(2)}|{\bf x}^{(1)})=\frac{f({\bf x})}{P({\bf X}^{(1)}={\bf x}^{(1)})},$$
the result follows from (1.4) and Theorem 4.1(i), since
$$|{\bf \Sigma}|=|{\bf \Sigma}_{11}|\cdot |{\bf \Sigma}_{22}-{\bf \Sigma}_{21}{\bf \Sigma}_{11}^{-1}{\bf \Sigma}_{12}|,$$
and  for each ${\bf x}\in\Bbb{R}^n$,
$${\bf x}'{\bf \Sigma}^{-1}{\bf x}={\bf x^{(1)}}'{\bf \Sigma}_{11}^{-1}{\bf x^{(1)}}+{\bf x_{2.1}}'{\bf \Sigma}_{22.1}^{-1}{\bf x_{2.1}}.$$
Here
$${\bf x_{2.1}}={\bf x^{(2)}}-{\bf \Sigma}_{21}{\bf \Sigma}_{11}^{-1} {\bf x^{(1)}}.$$
\begin{remark} Note that this conditional pdf is not a $(n-m)$-variate   multivariate elliptically symmetric logistic distribution unless
${\bf x^{(1)}}={\boldsymbol \mu}^{(1)}$.
\end{remark}

\section{Moments}
In this section we derive the moments of ${\bf X}$. From (1.7), we get, for real number $l>0$,
\begin{eqnarray*}
E(R^l)&=& \frac{1}{\int_0^{\infty}t^{\frac{n}{2}-1}g(t)dt}\int_0^{\infty}z^{\frac{n}{2}+l-1}g(z)dz\\
&=& \frac{1}{\int_0^{\infty}t^{\frac{n}{2}-1}g(t)dt}\int_0^{\infty}z^{\frac{n}{2}+l-1}\frac{e^{-bz}}{(1+e^{-az})^r} dz\\
&=& \frac{1}{\int_0^{\infty}t^{\frac{n}{2}-1}g(t)dt}\left(\frac{1}{a}\right)^{\frac{n}{2}+l}\int_0^{\infty}z^{\frac{n}{2}+l-1}\frac{e^{-\frac{b}{a}z}}{(1+e^{-z})^r}dz\\
&=& \frac{\left(\frac{1}{a}\right)^{l}\Gamma\left(\frac{n}{2}+l\right)\Phi_r^*\left(-1,\frac{n}{2}+l,\frac{b}{a}\right)}
{\Gamma\left(\frac{n}{2}\right)\Phi_r^*\left(-1,\frac{n}{2},\frac{b}{a}\right)},
\end{eqnarray*}
where    $\Phi^*_{r}$ is the generalized Hurwitz-Lerch Zeta function.

\begin{theorem}  Let ${\bf{X}}\sim GML_n ({\boldsymbol \mu},{\bf \Sigma},g)$, where $g$ is defined as (1.5).\\
(i) The expectation and the covariance are:
$$E({\bf X})={\boldsymbol \mu},\; Cov({\bf X})=\frac{\Phi_r^*\left(-1,\frac{n}{2}+1,\frac{b}{a}\right)}{2a\Phi_r^*\left(-1,\frac{n}{2},\frac{b}{a}\right)} {\bf \Sigma};$$
(ii)  For any integers $m_1,\cdots, m_n$, with $m=\sum_{i=1}^{n}m_i$, the product moments of   ${\bf Y}:={\bf \Sigma}^{-\frac12}({\bf X}-{\boldsymbol \mu})$ are
$$E\left(\prod_{i=1}^{n}Y_i^{2m_i}\right)= \frac{\left(\frac{1}{a}\right)^{m}\Gamma\left(\frac{n}{2}+m\right)\Phi_r^*\left(-1,\frac{n}{2}+m,\frac{b}{a}\right)}
{(\frac{n}{2})^{[m]}\Gamma\left(\frac{n}{2}\right)\Phi_r^*\left(-1,\frac{n}{2},\frac{b}{a}\right)}\prod_{i=1}^{n}\frac{(2m_i)!}{4^{m_i}(m_i)!},$$
 where $a^{[k]}=a(a+1)\cdots (a+k-1)$, $\Phi^*_{r}$ is the generalized Hurwitz-Lerch Zeta function.
\end{theorem}
{\bf Proof}  (i) By using (1.6) we have $E({\bf X})={\boldsymbol \mu}$ since $E{\bf U}^{(n)}=0$, and
$$ Cov({\bf X})=\frac{1}{n}E(R){\bf\Sigma}= \frac{\Phi_r^*\left(-1,\frac{n}{2}+1,\frac{b}{a}\right)}{2a\Phi_r^*\left(-1,\frac{n}{2},\frac{b}{a}\right)} {\bf \Sigma}.$$
(ii) By Eq. (2.18) in Fang et al. (1990), the product moments of $\bf Y$ are:
$$E\left(\prod_{i=1}^{n}Y_i^{2m_i}\right)=E(R^m)E\left(\prod_{i=1}^{n}U_i^{2m_i}\right).$$
The result follows since  (cf. Fang et al. (1990))
$$E\left(\prod_{i=1}^{n}U_i^{2m_i}\right)=\frac{1}{(\frac{n}{2})^{[m]}}\prod_{i=1}^{n}\frac{(2m_i)!}{4^{m_i}(m_i)!}.$$

\section{Characteristic function}

\begin{theorem}  Let ${\bf{X}}\sim GML_n ({\boldsymbol \mu},{\bf \Sigma},g)$, where $g$ is defined as (1.5).\\
The characteristic function of ${\bf X}$ can be expressed in the following  form:\\
(i) If $n=1$, assume that
 $$f(x)=\frac{d_1}{\sigma}g\left(\frac{(x-\mu)^2}{\sigma^2}\right), \;-\infty<x<\infty,$$
 where $\mu, \sigma>0$ are real number. Then
$$ \psi_{\bf X}(t)=d_1 e^{it\mu}\int_0^{\infty}\cos(t\sigma\sqrt{y})\frac{e^{-by}}{\sqrt{y}(1+e^{-ay})^r}dy,\; t\in(-\infty,\infty).$$
(ii) If $n\ge 2$, then
$$ \psi_{\bf X}({\bf t})=\frac{1}{ B(\frac{n-1}{2},\frac12)} \frac{1}{\Phi^*_{r}(-1,\frac{n}{2},\frac{b}{a})}e^{i{\bf t}'{\boldsymbol \mu}}\int_0^{\pi} \Phi^*_{r}\left(-1,\frac{n}{2},\frac{b}{a}-\frac{i}{a}{\bf t}'{\bf \Sigma}{\bf t}\cos\theta\right)\sin^{n-2}\theta  d\theta,$$
where $\Phi^*_{r}$ is the generalized Hurwitz-Lerch Zeta function.
\end{theorem}
{\bf Proof} (i) By definition, we have
\begin{eqnarray*}
\psi_{\bf X}(t)&=&\int_{-\infty}^{\infty}e^{itx}\frac{d_1}{\sigma}g\left(\frac{(x-\mu)^2}{\sigma^2}\right)dx\\
&=&2d_1e^{it\mu}\int_{0}^{\infty}\cos(t\sigma x)g(x^2)dx\\
&=&d_1e^{it\mu}\int_{0}^{\infty}\cos(t\sigma \sqrt{y})y^{-\frac12}g(y)dy\\
&=&d_1 e^{it\mu}\int_0^{\infty}\cos(t\sigma\sqrt{y})\frac{e^{-by}}{\sqrt{y}(1+e^{-ay})^r}dy.
\end{eqnarray*}
(ii) Using (1.6) and note that the independence of $R$ and  ${\bf U}^{(n)}$,
\begin{eqnarray*}
\psi_{\bf X}({\bf t})&=&E(e^{i{\bf t}'{\bf X}})=e^{i{\bf t}'{\boldsymbol \mu}} E(e^{i{\bf t}'{\sqrt{R}{\bf \Sigma}^{\frac12}{\bf U}^{(n)}}})\\
&=&e^{i{\bf t}'{\boldsymbol \mu}} E[E(e^{i({{\sqrt{R}{\bf \Sigma}^{\frac12} \bf t})'{\bf U}^{(n)}}}|R)]\\
&=&e^{i{\bf t}'{\boldsymbol \mu}} E[\Omega_n(R{\bf t}'{\bf \Sigma}{\bf t})]\\
&=&e^{i{\bf t}'{\boldsymbol \mu}}\int_0^{\infty} \Omega_n(v{\bf t}'{\bf \Sigma}{\bf t})P(R\in v)\\
&=&e^{i{\bf t}'{\boldsymbol \mu}}\int_0^{\infty} \Omega_n(v{\bf t}'{\bf \Sigma}{\bf t}) \frac{1}{\int_0^{\infty}t^{\frac{n}{2}-1}g(t)dt}v^{\frac{n}{2}-1}
\frac{\exp(-bv)}{(1+\exp(-av))^r}dv\\
&=& \frac{1}{\int_0^{\infty}t^{\frac{n}{2}-1}g(t)dt}e^{i{\bf t}'{\boldsymbol \mu}}\int_0^{\infty} \Omega_n(v{\bf t}'{\bf \Sigma}{\bf t})
\frac{v^{\frac{n}{2}-1}\exp(-bv)}{(1+\exp(-av))^r}dv\\
&=&\frac{1}{\int_0^{\infty}t^{\frac{n}{2}-1}g(t)dt}e^{i{\bf t}'{\boldsymbol \mu}}\int_0^{\infty}\left(\frac{1}{B(\frac{n-1}{2},\frac12)}\int_0^{\pi}
\exp(iv{\bf t}'{\bf \Sigma}{\bf t}\cos\theta)\sin^{n-2}\theta d\theta\right)
\frac{v^{\frac{n}{2}-1}\exp(-bv)}{(1+\exp(-av))^r}dv\\
&=& \frac{1}{B(\frac{n-1}{2},\frac12)} \frac{1}{\int_0^{\infty}t^{\frac{n}{2}-1}g(t)dt}e^{i{\bf t}'{\boldsymbol \mu}}\int_0^{\pi} \sin^{n-2}\theta\left(\int_0^{\infty}
\exp(iv{\bf t}'{\bf \Sigma}{\bf t}\cos\theta)
\frac{v^{\frac{n}{2}-1}\exp(-bv)}{(1+\exp(-av))^r}dv\right) d\theta\\
&=& \frac{1}{ B(\frac{n-1}{2},\frac12)} \frac{1}{\Phi^*_{r}(-1,\frac{n}{2},\frac{b}{a})}e^{i{\bf t}'{\boldsymbol \mu}}\int_0^{\pi} \Phi^*_{r}\left(-1,\frac{n}{2},\frac{b}{a}-\frac{i}{a}{\bf t}'{\bf \Sigma}{\bf t}\cos\theta\right)\sin^{n-2}\theta  d\theta,
\end{eqnarray*}
where $\Omega_n({\bf t}'{\bf t})$ is the characteristic function of ${\bf U}^{(n)}$ (see Fang et al. (1990, (3.1))):
 $$\Omega_n(||{\bf t}||^2)=\frac{1}{B(\frac{n-1}{2},\frac12)}\int_0^{\pi}
\exp(i||{\bf t}||\cos\theta)\sin^{n-2}\theta d\theta.$$
\begin{remark}
Using the following equivalent forms of $\Omega_n(||{\bf t}||^2)$ (see  Fang et al. (1990, (3.2), (3.3))
$$\Omega_n(||{\bf t}||^2)=\frac{\Gamma(\frac{n}{2})}{\sqrt{\pi}\Gamma(\frac{n-1}{2})}\sum_{k=0}^{\infty}\frac{(-1)^k||{\bf t}||^{2k}\Gamma(\frac{n-1}{2})
\Gamma(\frac{2k+1}{2})}{(2k)!\Gamma(\frac{n+2k}{2})}
$$
and
$$\Omega_n(||{\bf t}||^2)={}_0F_1\left(\frac{n}{2};\frac14||{\bf t}||^{2}||\right),$$
 we get
the following equivalent forms of   the characteristic function  $\psi_{\bf X}({\bf t})$:
$$\psi_{\bf X}({\bf t})=\frac{e^{i{\bf t}'{\boldsymbol \mu}}\Gamma(\frac{n}{2})}{\sqrt{\pi}\Phi^*_{r}\left(-1,\frac{n}{2},\frac{b}{a}\right)}\sum_{k=0}^{\infty}\frac{(-1)^k({\bf t}'{\bf \Sigma}{\bf t} )^k \Gamma(\frac{2k+1}{2})}{(2k)!\Gamma(\frac{n+2k}{2})}\Phi^*_{r}\left(-1,\frac{n}{2}+k,\frac{b}{a}\right),$$
and
 $$\psi_{\bf X}({\bf t})=\frac{e^{i{\bf t}'{\boldsymbol \mu}}}{\Phi^*_{r}\left(-1,\frac{n}{2},\frac{b}{a}\right)}\sum_{k=0}^{\infty}\frac{(\frac14{\bf t}'{\bf \Sigma}{\bf t} )^k}{(\frac{n}{2})^{[k]}k!}\Phi^*_{r}\left(-1,\frac{n}{2}+k,\frac{b}{a}\right),$$
where ${}_0F_1(\cdot;\cdot)$  is the  generalized hypergeometric function, $\Phi^*_{r}$ is the generalized Hurwitz-Lerch Zeta function and $a^{[k]}$  are the ascending factorials, i.e. $a^{[k]}=a(a+1)\cdots (a+k-1).$
\end{remark}

\noindent{\bf\Large Acknowledgements} \ 

\noindent The research was supported by the National Natural Science Foundation of China (No. 11171179, 11571198).

\end{document}